\documentclass{amsart}
\usepackage{amssymb,amsmath}
\newcommand{\RR}{{\mathbb R}}

\newcommand{\e}{\varepsilon}

\newcommand{\om}{\omega}

\newcommand{{\loc}}{{\ell\mathrm oc}}

\def\meanint{{\diagup\hskip -.42cm\int}}

\newtheorem{theorem}{Theorem}
\newtheorem{lemma}{Lemma}
\newtheorem{proposition}{Proposition}
\newtheorem{corollary}{Corollary}

\newtheorem{remark}{Remark}

\begin{document}

\title{ Gilbarg-Serrin Equation and Lipschitz Regularity}
\author{Vladimir Maz'ya}
\address{Department of Mathematics, Link\"oping University, SE-581 83 Link\"oping, Sweden and RUDN, 6 Miklukho-Maklay St, Moscow, 117198, Russia}
\email{vladimir.mazya@liu.se}

\author{ Robert McOwen}
\address{Department of Mathematics, Northeastern University, Boston, MA 02115}
\email{r.mcowen@northeastern.edu}
\date{September 11, 2021}
\dedicatory{Dedicated in memory of James Serrin who passed away ten years ago.}


\begin{abstract}
We discuss conditions for Lipschitz and $C^1$ regularity of solutions for a uniformly elliptic equation in divergence form. We focus on  coefficients having the form that was introduced by Gilbarg \& Serrin. In particular, we find cases where Lipschitz or $C^1$ regularity holds but the coefficients are not Dini continuous, or do not even have Dini mean oscillation. The form of the coefficients also enables us to obtain specific conditions and examples for which there exists a weak solution that is not Lipschitz continuous.

\end{abstract}
\maketitle

\section{Introduction and General Theory}\label{sec:Intro}

The topic of this paper is Lipschitz and $C^1$ regularity for weak solutions of a uniformly elliptic equation in divergence form:
\begin{equation}\label{Lu=0}
\partial_j (a_{ij}(x)\,\partial_i u)=0, \quad\hbox{for $x\in \Omega$},
 \end{equation}
 where the coefficients $a_{ij}=a_{ji}$ are bounded, measurable functions in an open set $\Omega\subset\RR^n$, $n\geq 2$, which contains the origin $x=0$.  For a restricted class of coefficients $a_{ij}$, we
seek conditions that guarantee that any weak solution $u\in H_{\loc}^{1,2}(\Omega)$, i.e.\ 
 $u$ and $\nabla u$ are locally in $L^2$, must be {\it Lipschitz continuous at $x=0$},
 i.e.\ $|u(x)-u(0)|\leq c\,|x|$ holds for $x\in B_\rho$ with $\rho$ sufficiently small,  or {\it Lipschitz continuous in a neighborhood of $x=0$}, i.e.\ $|u(x)-u(y)|\leq c\,|x-y|$ for  $x,y\in B_\rho$ with $\rho$ sufficiently small.
 (Here and throughout this paper, $c$ denotes a positive constant and $B_\rho$ denotes the ball of radius $\rho$ centered at $x=0$ such that $\overline{B_\rho}\subset\Omega$; the values of $c$ and $\rho$ may change with each occurrence.)
 Recall that Lipschitz continuous functions are differentiable almost everywhere (by Rademacher's theorem) and so our solution will satisfy $\nabla u\in L^\infty_{\loc}(\Omega)$. We also seek conditions implying that any weak solution is differentiable at $x=0$, or
 is $C^1$ near $x=0$.

\subsection{Background}
Let us review a little of what is known about the regularity of weak solutions of \eqref{Lu=0}. 
Any weak solution of \eqref{Lu=0} is known to be H\"older continuous, i.e.\  $|u(x)-u(y)|\leq c\,|x-y|^\alpha$ for some $\alpha\in (0,1)$, at least locally in $\Omega$ (cf.\ De Giorgi \cite{DG}, Nash \cite{N}, Moser \cite{Mo}, Landis \cite{L}).
 When the coefficients  are continuous in $\Omega$, then it is well-known (cf.\ Agmon, Douglis, Nirenberg \cite{ADN}) that $\nabla u\in L^p_{\loc}(\Omega)$ for $1<p<\infty$; in fact, this is even true when the coefficients are in VMO (cf.\ Chiarenza, Frasca, Longo \cite{CFL}, Di Fazio \cite{DF}), or in BMO but sufficiently close (depending on $pp'$) to VMO  (cf.\ Maz'ya, Mitrea, Shaposhnikova \cite{MMS}).\footnote{Here, as usual, BMO stands for bounded mean oscillation and VMO stands for vanishing mean oscillation.}
Note that continuity of the coefficients is not sufficient to show $\nabla u\in L^\infty_{\loc}(\Omega)$ (cf.\ Jin,  Maz'ya, Van Schaftingen \cite{JMV}).
 However, if the coefficients are Dini continuous in $\Omega$, then a weak solution is known to be not only Lipschitz continuous but  $C^1$ (cf.\  Hartman,  Wintner \cite{HW}, Burch \cite{B}, Taylor \cite{T}).
Special attention has been given to pointwise bounds on the gradient of the solutions for linear and nonlinear elliptic equations and systems in \cite{MM}, \cite{CM}, \cite{KuM}, \cite{CM1}, \cite{KuM2}, \cite{CM2}, \cite{CM3}. The Dini condition also provides regularity in parabolic problems: \cite{ME}, \cite{I}, \cite{A}.
 
  The classical Dini condition on the coefficients involves the modulus of continuity $\omega$:
 \begin{equation}\label{Dini1}
 |a_{ij}(x)-a_{ij}(y)|\leq  \omega(|x-y|) \quad \hbox{for $x,y\in\Omega$,}
 \end{equation}
 where $\omega(r)$ is a continuous, nondecreasing function $[0,\infty)\to [0,\infty)$ with $\om(0)=0$. The coefficients are said to be {\it Dini continuous} if $\omega$ satisfies
 \begin{equation}\label{Dini2}
 \int_0^{1/2} \frac{\omega(r)}{r}\,dr <\infty.
 \end{equation}
 When $\omega(r)=r^\alpha$ for some $\alpha\in (0,1)$ then this regularity coincides with H\"older continuity. There are weaker moduli of continuity that satisfy the Dini condition. If we let $\omega(r)=|\log r|^{-\gamma}$ for $0<r<1/2$ and $\omega(0)=0$, then $\omega$ satisfies \eqref{Dini2} if  $\gamma>1$. However,
 if $\omega(r)=|\log r|^{-\gamma}$ for $0<\gamma\leq 1$, then the coefficients are not Dini continuous, and the usual regularity properties of solutions may not hold.

More recently, generalizations of the Dini condition have been used. For example, consider the ``square-Dini condition'' on the modulus of continuity $\omega$:
  \begin{equation}\label{DiniSquared}
 \int_0^{1/2} \frac{\omega^2(r)}{r}\,dr <\infty.
 \end{equation}
 A function is {\it square-Dini continuous} if its modulus of continuity $\omega$ satisfies \eqref{DiniSquared}.
 Clearly the condition \eqref{DiniSquared}  is more general than \eqref{Dini2}; for 
 example, $\omega(r)=|\log r|^{-\gamma}$ satisfies  \eqref{DiniSquared}  when $\gamma>1/2$.
 In  \cite{MM}, the square-Dini condition was only required at a point (see Section 1.2) and was used to
 reduce the regularity of a weak solution at that point to stability properties of a dynamical system obtained from the coefficients. 
This dynamical systems approach to studying the properties of solutions of elliptic equations has been used in many other publications: \cite{KM}, \cite{KM2}, \cite{Mc}, \cite{MM1}, \cite{MM2}. We shall have more to say about this dynamical system method below, but for now let us mention another generalization of the Dini condition.
 
 Several authors have  considered a ``Dini mean oscillation'' condition on the coefficients.
 In the comprehensive paper \cite{KuM},
 Kuusi \& Mingione used the notion of Dini mean oscillation in their unified linear/nonlinear potential theory to obtain gradient bounds for solutions of quasilinear elliptic equations.
For linear equations the condition was used by Dong \& Kim \cite{DK} to show that weak solutions of an elliptic equation in divergence form are $C^1$ and for nondivergence form they are $C^2$. Let us formulate the condition on the coefficients
as follows: a function $a$ that is integrable in $\Omega$ has {\it Dini mean oscillation} if the following function $\om_a$ satisfies the Dini condition \eqref{Dini2}:
 \begin{equation}
 \om_a(r):=\sup_{B_r(x)\subset\Omega} \ \meanint_{B_r(x)} |a(y)-\widetilde a_{x,r}|\,dy
 \quad
 \hbox{where} \ \widetilde a_{x,r}:= \meanint_{B_r(x)} a(y)\,dy.
 \end{equation}
 (Here $B_r(x)=\{y: |x-y|<r\}$
 and the slashed integral denotes mean value.)
 There are functions that are not Dini continuous but do have Dini mean oscillation; for example, if $a(x)=(-\log|x|)^{-\gamma}$ in
 $\Omega=B_{1/2}$, then $\omega_a(r)=c\,|\log r|^{-\gamma-1}$, so $a$ has Dini mean oscillation for all $\gamma>0$. The paper \cite{DK} has spawned a series of applications to various linear elliptic and parabolic equations: \cite{DEK1}, \cite{DEK2}, \cite{DLK}. As useful as the concept of Dini mean oscillation has proven to be, however, we will see that there are equations
of the form \eqref{Lu=0} with coefficients that do not have Dini mean oscillation, yet Lipschitz and $C^1$ regularity can be obtained by other methods.
 
 \subsection{Dynamical Systems Method}
 Now let us describe the method in \cite{MM} of using dynamical systems to study the regularity  at a given point
 of weak solutions of \eqref{Lu=0}. We remark that the conditions on the coefficients need only occur at the given point. These
conditions and the dynamical system are most conveniently formulated when the point in question is $x=0$ and $a_{ij}(0)=\delta_{ij}$; this can always be achieved by a change of coordinates at $0$ since we have assumed the coefficients are symmetric. 
 Under these assumptions, the condition on the coefficients is
 \begin{equation}\label{ptwiseDini^2}
 \sup_{|x|=r}|a_{ij}(x)-\delta_{ij}|\leq \om(r) \quad\hbox{as}\ r\to 0,
 \end{equation}
 where $\om$ satisfies the square-Dini condition \eqref{DiniSquared}; we describe this by saying the coefficients are 
 {\it square-Dini continuous at $x=0$.}
 
 The basic idea  in \cite{MM} is to write the solution as 
 \begin{equation}\label{eq:u(x)=}
 u(x)=u_0(|x|)+\vec v(|x|)\cdot x+w(x),
 \end{equation}
 where
 \[
 u_0(r)=\meanint_{S^{n-1}}u(r\theta)\,d\theta \quad\hbox{and}\quad \vec v_k(r)=\frac{n}{r}\meanint_{S^{n-1}} u(r\theta)\theta_k\,d\theta;
 \]
 here $r=|x|$, $\theta=x/|x|\in S^{n-1}$, and $d\theta$ is standard surface measure on $S^{n-1}$. 
 As shown in  \cite{MM}, if the coefficients are square-Dini continuous at $x=0$, then the behaviors of both $\vec v(r)$ and $r\,\vec v\,'(r)$ are controlled by the asymptotic properties of the solutions of the $(n\times n)$ dynamical system 
 \begin{equation}
 \frac{d\phi}{dt} + R(e^{-t})\,\phi=0 \quad\hbox{for}\ T<t<\infty,
 \label{eq:DynSystm}
 \end{equation}
where $t=-\log r$, $T$ is sufficiently large, and $R$ is the ($n\times n$)-matrix
 \begin{equation}
 R(r):=\meanint_{S^{n-1}}\left(A(r\theta)-nA(r\theta)\theta\otimes\theta\right)\,d\theta, \quad\hbox{for $0<r<1$,}
\label{def:R}
\end{equation}
with $A=(a_{ij})$ and $A\theta\otimes\theta$ the outer product of the vectors $A\theta$ and $\theta$. 
Note that  $\pmb{|}R(r)\pmb{|}\leq c\,\om(r)$, where  $\pmb{|}\cdot\pmb{|}$ denotes the matrix norm; also note that $R$ need not be symmetric. The asymptotic property that we require for solutions to be Lipschitz is {uniform stability}: we say that 
\eqref{eq:DynSystm} is {\it uniformly stable as $t\to\infty$} if for every $\e>0$ there exists $\delta=\delta(\e)>0$ such that any solution of \eqref{eq:DynSystm} satisfying $|\phi(t_1)|<\delta$ for some $t_1>T$ satisfies $|\phi(t)|<\e$ for all $t>t_1$. We show in \cite{MM} that the uniform stability of \eqref{eq:DynSystm} implies not only that  $\vec v(r)$ and $r\,\vec v\,'(r)$ are bounded as $r\to 0$, but with additional analysis that
$|u_0(r)-u_0(0)|,|w(x)|=O(r\,\om(r))$ as $r\to 0$. This means that
\begin{equation}\label{eq:u(x)=}
 u(x)=u_0(0)+\vec v(|x|)\cdot x+O(r\,\om(r)),
 \end{equation}
 which confirms that $u$ is Lipschitz at $x=0$. 

For differentiabilty, we also require that every solution of \eqref{eq:DynSystm} be {\it asymptotically constant}, i.e.\ $\phi(t)\to\phi_\infty$ as $t\to\infty$. In that case, we have $\vec v(r)=\vec v(0)+o(1)$ as $r\to 0$, which shows that $u$ is differentiable at $x=0$.
Here is the complete result: 
\begin{theorem}[Theorem 1 in \cite{MM}]\label{th:1} 
Suppose the $a_{ij}$ are bounded, measurable functions in $\Omega$ that are square-Dini continuous at $x=0$. Suppose also  that the dynamical system \eqref{eq:DynSystm} is uniformly stable. Then every weak solution $u\in H^{1,2}(\Omega)$ of \eqref{Lu=0} is Lipschitz continuous at $x=0$ and satisfies
\begin{equation}\label{est:u(x)-u(0)}
|u(x)-u(0)|\leq \frac{c\,|x|}{\rho}\left(\meanint_{|y|<\rho} |u(y)|^2\,dy\right)^{1/2} \quad\hbox{for}\ |x|<\rho/2,
\end{equation}
 for $\rho$ sufficiently small. In addition, if every solution of \eqref{eq:DynSystm}  is {asymptotically constant}, then $u$ is differentiable at $x=0$ and
 \begin{equation}\label{nabla_u(x)=}
\partial_j u(0)=\lim_{r\to 0} \frac{n}{r} \,\meanint_{S^{n-1}} u(r\theta)\theta_j d\theta, \quad \hbox{where}\ \theta_j=x_j/r.
\end{equation}
\end{theorem}

We also gave in \cite{MM} several simpler sufficient conditions for the uniform stability of \eqref{eq:DynSystm}, such as 
\begin{equation}
 \frac{R(r)}{r} \int_0^r \frac{R(\rho)}{\rho}\,d\rho \in L^1(0,\e),
 \end{equation}
or 
\begin{equation}
\begin{aligned}
&\int_{r_1}^{r_2}\frac{\mu({\mathcal S}(\rho))}{\rho}\,d\rho \leq K<\infty \quad\hbox{for all}\ 0<r_1<r_2<\e, \\
&  \hbox{ where $\mu({\mathcal S})=$ largest eigenvalue of ${\mathcal S}:=-(R+R^t)/2$.} 
\end{aligned}
 \end{equation}
 The condition that all solutions of a dynamical system are asymptotically constant is independent of its uniform stability (cf.\ Section 5 in \cite{MM}); when \eqref{eq:DynSystm} is a scalar equation, as we shall have in Section 2.1, it is elementary to completely characterize this condition.
  
Also note that, if the $a_{ij}$ are radial functions, then $R= 0$ and Theorem \ref{th:1} implies Lipschitz regularity and differentiability at $x=0$ only assuming the $a_{ij}$ are square-Dini continuous at $x=0$. More generally, if we let $a^0_{ij}(r)$ denote the mean of $a_{ij}$ over the sphere $|x|=r$, then we can uniquely write  $a_{ij}(x)=a^0_{ij}(|x|)+a^1_{ij}(x)$; if $a^1_{ij}(x)$ is Dini continuous and $a_{ij}^0(|x|)$ is only square-Dini continuous at $x=0$, then every solution is Lipschitz continuous and differentiable at $x=0$.

\smallskip
 \subsection{From a Point to a Neighborhood}
Theorem \ref{th:1} above concerns the Lipschitz continuity and differentiability at $x=0$, but under certain conditions, these 
regularity properties extend locally.
Let us consider a simple case where the coefficients  are almost everywhere differentiable with an estimate for the growth of $\nabla a_{ij}(x)$ as $|x|\to 0$ and obtain a uniform bound near $x=0$ on the gradient of a weak solution; this will also show the solution is $C^1$ near $x=0$.
\begin{theorem}\label{th:2}  
Suppose the $a_{ij}$ are differentiable in $\Omega\backslash\{0\}$ with 
\begin{equation}\label{est:grad-a}
|\nabla a_{ij}(x)| \leq \frac{c}{|x|} \quad\hbox{for $0<|x| <2\rho$.}
\end{equation}
Suppose the $a_{ij}$ are square-Dini continuous at $x=0$
and the dynamical system  \eqref{eq:DynSystm} is uniformly stable.
Then every weak solution $u\in H^{1,2}(\Omega)$ of \eqref{Lu=0}
is Lipschitz continuous in $B_\rho$ and there is a positive constant $c$ such that
\begin{equation} \label{est:grad-u,rho}
|\nabla u(x)|\leq \frac{c}{\rho}\left(\meanint_{|y|<2\rho} |u(y)|^2\,dy\right)^{1/2} \quad\hbox{for $0<|x|<\rho$.}
\end{equation}
Moreover, $u\in C^1(\overline{B_\rho})$.
\end{theorem}

\noindent{\bf Proof}.
We will require standard elliptic interior estimates in an annulus. For $\rho>0$, let us introduce the annulus
\begin{equation}
A_{r}=\{x\in\RR^n: \frac{r}{2}<|x|<2r\}.
\end{equation}
We first consider $r=1$. We know $|\nabla a_{ij}(x)|\leq c$ for $x\in A_{1}$, so if
$A_1'$ is a concentric open sub-annulus of $A_1$ such that $\overline {A_1'}\subset A_{1}$ and  $1<p<\infty$,
then there is a positive constant $c$ such that
\begin{equation} \label{est:H2p-u}
\| u \|_{H^{2,p}(A_1')} \leq c\,\|u\|_{L^1(A_{1})}.
\end{equation}
(The above can be obtained from the standard $L^p$ elliptic estimates: cf.\ Theorem 15.1'' 
and Remark 1 following Theorem 7.3 in \cite{ADN};
also (14.3.6) in Chapter 14 of \cite{MS}.)
Taking $p>n$ so that we can use the Sobolev embedding $H^{2,p}(A_1')\subset C^1(\overline{A_1'})$, we obtain the following:
\begin{equation} \label{est:sup-grad-u}
|\nabla u(x) |\leq c\,\|u\|_{L^1(A_{1})} \quad\hbox{for}\ \frac{3}{4}\leq |x|\leq \frac{3}{2}.
\end{equation}
Applying a dilation to the annulus, we find that if
\begin{equation} \label{est:sup-grad-aij-r}
|\nabla a_{ij}(x)| \leq \frac{c}{r} \quad\hbox{for $x\in A_r$,}
\end{equation}
where $c$ is independent of $r$, then  there is a constant $c$ independent of $r$ such that
\begin{equation} \label{est:sup-grad-u-r}
|\nabla u(x)|\leq \frac{c}{r} \,\meanint_{A_r}|u(y)|\,dy
 \quad\hbox{for $\frac{3}{4}r \leq |x| \leq \frac{3}{2}r$.}
\end{equation}

For $x\in B_\rho\backslash\{0\}$, choose $r>0$ so that $\frac{3}{4} r\leq |x|\leq \frac{3}{2}r$.
Since $u(x)-u(0)$ is also a solution of \eqref{Lu=0}, we can apply \eqref{est:sup-grad-u-r} to it and then use
\eqref{est:u(x)-u(0)}:
\begin{equation} \label{est:sup-grad-u-rho}
\begin{aligned}
|\nabla u(x)|&\leq \frac{c}{r}\,\meanint_{A_r}|u(y)-u(0)|\,dy 
 \leq \frac{c}{r}\,\meanint_{A_r}\frac{c|y|}{\rho} dy  \left(\meanint_{|z|<2\rho}|u(z)|^2\,dz\right)^{1/2} \\
&\leq \frac{c}{\rho} \left(\meanint_{|z|<2\rho}|u(z)|^2\,dz\right)^{1/2},
\end{aligned}
\end{equation}
where we have used $|y|<2r$ for $y\in A_r$. 

The uniform bound \eqref{est:grad-u,rho} also enables us to conclude that $u\in C^1(\overline{B_\rho})$. In fact, we can use a mollifier  to obtain coefficients $a_{ij}^h(x)$ that are smooth on $B_{2\rho}$ and $a_{ij}^h\to a_{ij}$ in $L^p(B_{2\rho})$ as $h\to 0$ for $1\leq p<\infty$.  If we let $u^h$ denote the solution of \eqref{Lu=0} in $B_{2\rho}$ with coefficients $a_{ij}^h$ and Dirichlet condition $u^h=u$ on $\partial B_{2\rho}$, then we know that $u^h\in C^\infty(B_{2\rho})$ and $u^h\to u$ in $H^{1,2}(B_{2\rho})$  (see pages 150-151 in  \cite{M} or Section 4 of Chapter 2 in \cite{L}).
Moreover the continuity modulus of the $a_{ij}^h$ does not increase, and the estimate  \eqref{est:grad-a} continues to hold with a constant $c$ independent of $h$ near zero, so  the estimate \eqref{est:grad-u,rho} applies to $u^h$ with $c_1$ independent of $h$  near zero. If we consider a sequence $v^m=u^{h_m}$ with $h_m\to 0$ as $m\to\infty$, we know  $\|v^{\ell}-v^{m}\|_{L^2(B_{2\rho})}\to 0$ as $\ell,m\to\infty$  and can use \eqref{est:grad-u,rho} to conclude  $\|\nabla(v^{\ell}-v^{m})\|_{L^\infty(B_\rho)}\to 0$. By completeness of the Banach space $C^1(\overline{B_\rho})$, this shows $u\in C^1(\overline{B_\rho})$.
$\Box$

\begin{remark}\label{re:1} 
Instead of assuming differentiability of the coefficients for $x\not=0$ in Theorem 1, we could have assumed they
satisfy $|a_{ij}(x)-a_{ij}(y)|\leq \omega(|x-y|/|x|)$ for $|x-y|<\frac{1}{2}|x|$, where $\omega$ satisfies the Dini condition \eqref{Dini1}.
Let us also mention the recent work of De Filippis \& Mingione \cite{DM} where \eqref{est:grad-a} and our additional conditions on the coefficients at $x=0$ are  replaced by
the gradient of the coefficients belonging to the Lorentz space $L(n,1)$.
\end{remark}

 \section{Gilbarg-Serrin Equations}\label{sec:GS}
 
 In \cite{GS} and \cite{Se}, Gilbarg and Serrin considered coefficients of the form
 \begin{equation}\label{GS}
 a_{ij}(r,\theta)=\delta_{ij}+g(r)\,\theta_i\theta_j, \quad\hbox{where $r=|x|$ and $\theta_i=x_i/r$.}
 \end{equation}
They assumed $g(r)$ is a bounded function  satisfying $g(r)>-1+\e$, which guarantees uniform ellipticity,
and studied solutions of  equations involving the associated operators in both divergence and nondivergence form.
In this section, we consider conditions on $g(r)$ that guarantee that any weak solution of \eqref{Lu=0} must be Lipschitz continuous or differentiable at $x=0$,  or when there exists a weak solution that is not Lipschitz continuous at $x=0$.
In cases where the solution is Lipschitz continuous in a neighborhood of $x=0$, we will use Theorem \ref{th:2} to conclude that the solution is $C^1$.
Let us state several results in this direction. 

\subsection{Application of the Dynamical Systems Method}
Let us apply the results of the dynamical systems method to coefficients of the form \eqref{GS}.
It is shown in Section 5 of \cite{MM} that the dynamical system \eqref{eq:DynSystm}  reduces in this case to the (scalar) ordinary differential equation
\begin{equation}
\frac{d\phi}{dt}=\frac{n-1}{n}\, {\rm g}(t)\,\phi,
\label{eq:ODE}
 \end{equation}
where ${\rm g}(t)=g(e^{-t})$. We can integrate this ODE to see that it is uniformly stable if and only if
there is a real constant $K$ such that
\[
\int_s^t {\rm g}(\tau)\,d\tau \leq K <\infty \quad\hbox{for $t>s>T$,}
 \]
and all solutions are asymptotically constant if and only if the improper integral
\[
\int_T^\infty {\rm g}(\tau)\,d\tau \quad \hbox{converges to an extended real number $<\infty$.}
 \]
We are now able to combine Theorems \ref{th:1} and \ref{th:2} to obtain:

\begin{proposition}\label{pr:1} 
Let the coefficients $a_{ij}$ be given by  \eqref{GS}, where $g(r)$ is square-Dini continuous at $0$, and let $u\in H^{1,2}(\Omega)$ be a weak solution of \eqref{Lu=0}. If, for some real constant $K$, we have
\begin{equation}\label{cond:g(r)/r}
\int_{r_1}^{r_2} \frac{g(r)}{r}\,dr \leq K<\infty \quad\hbox{for \ $0<r_1<r_2<\rho$,}
\end{equation}
then $u$ is Lipschitz continuous at $x=0$ and satisfies \eqref{est:u(x)-u(0)}.
If, in addition, 
\begin{equation}\label{cond:asym_const}
\int_0^1  \frac{g(r)}{r}\,dr  \quad \hbox{converges to an extended real number $<\infty$,}
\end{equation}
then $u$ is differentiable at $x=0$ and satisfies \eqref{nabla_u(x)=}.
On the other hand, if we just assume \eqref{cond:g(r)/r} and that $g(r)$ is differentiable for $r>0$ with derivative satisfying $|r\,g'(r)|\leq c$ for $0<r<\rho$, then $u\in C^1(\overline{B_\rho})$ and satisfies the estimate \eqref{est:grad-u,rho}.
\end{proposition}

\begin{corollary}
If $g(r)$ is square-Dini continuous at $r=0$ and satisfies both \eqref{cond:g(r)/r} and 
\begin{equation}
\int_\e^1\frac{g(r)}{r}\,dr\to -\infty\quad\hbox{as}\ \e\to 0,
\end{equation}
then every weak solution of \eqref{Lu=0} is differentiable at $r=0$ with $\nabla u(0)=0$.
\end{corollary}

\subsection{Application of a Comparison Principle}
We will discuss a comparison principle based upon solutions of the following ordinary differential equation:
\begin{equation}\label{ode:v}
\frac{1}{r^{n-1}}\frac{d}{dr}\left( r^{n-1} (1+g(r))\frac{dv}{dr}\right)-\frac{n-1}{r^2}\,v=0, \quad\hbox{for\ $0<r<1$.}
\end{equation}
\begin{proposition}\label{pr:2} 
Let $Z(r)$ denote a solution of \eqref{ode:v} satisfying the finite energy condition
\begin{equation}\label{finite_energy}
\int_0^1 \left( |Z'(r)|^2+r^{-2}|Z(r)|^2\right)r^{n-1}\,dr <\infty.
\end{equation}
If $Z(r)$ satisfies $0<Z(r)\leq c\,r$ as $r\to 0$, then every weak solution $u\in H^{1,2}(\Omega)$ of \eqref{Lu=0} with coefficients given by \eqref{GS} is Lipschitz continuous at $x=0$. In addition, if $g(r)$ in \eqref{GS} is differentiable for $r>0$ and 
$|rg'(r)|\leq c$ for $0<r<2\rho$, then $u\in C^1(\overline{B_\rho})$ and satisfies the estimate \eqref{est:grad-u,rho}.

\end{proposition}

\noindent
{\bf Proof.} We let $Z(r)$ be a solution of  \eqref{ode:v} satisfying $Z(1)>0$ and the finite energy condition \eqref{finite_energy}. By the non-oscillatory character of \eqref{ode:v}, $Z(r)$ cannot change sign for $0<r<1$, so we have $Z(r)>0$. With coefficients $a_{ij}$ given by \eqref{GS}, the equation  \eqref{Lu=0} can be written as
\begin{equation}\label{GS:Lu=0}
\frac{1}{r^{n-1}}\frac{\partial}{\partial r} \left(r^{n-1}(1+g(r)) \frac{\partial u}{\partial r}\right)+\frac{1}{r^2}\Delta_\theta u=0,
\end{equation}
where $\Delta_\theta$ denotes the Laplace-Beltrami operator on the sphere $S^{n-1}$. If  $u\in H^{1,2}(B_1)$ is a weak solution, we can integrate this equation over $S^{n-1}$ to find that $u$ satisfies
\[
\meanint_{S^{n-1}} u(r\,\theta)\,d\theta= constant = u(0),
\]
where $u(0)$ denotes the value of $u(x)$ at $x=0$. But since we can subtract a constant from a weak solution of \eqref{Lu=0}, we may assume that $u(0)=0$, so $u(r,\theta)$ is orthogonal to 1 on every sphere centered at $0$:
\begin{equation}\label{u:orthog}
\int_{S^{n-1}} u(r\,\theta)\,d\theta=0 \quad\hbox{for $0<r<1$},
\end{equation}
We use the following comparison principle between $Z$ and $u$: 
\begin{lemma}\label{le:comparison} 
Let the coefficients $a_{ij}$ be given by \eqref{GS} and  $u\in H^{1,2}(B_1)$ is a weak solution of \eqref{Lu=0} satisfying \eqref{u:orthog}.  Then
\begin{equation}\label{eq:comparison}
\left( \meanint_{S^{n-1}} |u(r\,\theta)|^2\,d\theta\right)^{1/2} \leq \frac{Z(r)}{Z(1)}  \left(\meanint_{S^{n-1}} |u(\theta)|^2\,d\theta\right)^{1/2}
\quad\hbox{for\ $0<r<1$.}
\end{equation}
\end{lemma}
\noindent
We will prove Lemma \ref{le:comparison} below. But to complete the proof of Proposition \ref{pr:2},
recall  from Moser \cite{Mo} that a weak solution of  \eqref{Lu=0} satisfies the estimate
\begin{equation}
\sup_{|x|<\delta} |u(x)|\leq c\left(\meanint_{B_{2\delta}} |u(y)|^2\,dy\right)^{1/2} \quad\hbox{for}\ 0<\delta\leq 1/2.
\end{equation}
We may combine this with \eqref{eq:comparison} and the assumption $Z(r)\leq c\,r$ to conclude that $u(x)$ is Lipschitz continuous at $x=0$. If $g$ is differentiable with $|rg'(r)|\leq c$ for $0<r<2\rho$, then we can use Theorem \ref{th:2} to conclude
$u\in C^1(\overline{B_\rho})$ and obtain the estimate \eqref{est:grad-u,rho}.
$\Box$

\begin{remark}\label{re:1.5} 
The differentiability at $x=0$ of the solution $Z(r)\Theta_1(\theta)$ of \eqref{Lu=0} means
$Z(r)\Theta_1(\theta)=V\cdot x+o(r)$ where $V$ is a constant vector. Multiplying both sides by $\Theta_1$ and integrating over the unit sphere
we obtain $Z(r)=c\,r+o(r)$ for some constant $c$.
\end{remark}

\medskip \noindent  
{\bf Proof of Lemma \ref{le:comparison}.} Considering $u$ as a weak solution of \eqref{GS:Lu=0} in the ball $B_\rho$ for $0<\rho<1$, we have
\begin{equation}
\int_{B_\rho(0)} \left((1+g(r))\left|\frac{\partial u}{\partial r}\right|^2+\frac{1}{r^2}\left|\nabla_\theta u\right|^2\right)\,dx
=(1+g(\rho))\rho^{n-1}\int_{S^{n-1}} u(\rho\,\theta)\frac{\partial u}{\partial r}(\rho\,\theta)\, d\theta.
\end{equation}
Using \eqref{u:orthog}, minimizing the integration over the sphere $S^{n-1}$ on the left hand side, and rewriting the expression on the right hand side, we find
\begin{equation}\label{estimate}
\int_{S^{n-1}}\! \int_0^\rho \left((1+g(r))\left|\frac{\partial u}{\partial r}\right|^2+\frac{n-1}{r^2}\left| u\right|^2\right)r^{n-1}drd\theta
\leq (1+g(\rho))\frac{\rho^{n-1}}{2 }\frac{d}{d \rho}\int_{S^{n-1}}(u(\rho\,\theta))^2d\theta.
\end{equation}
Now consider the variational functional associated with \eqref{ode:v}:
\begin{subequations}\label{variational}
\begin{equation}\label{eq:J=}
J(v)=\int_0^\rho \left( (1+g(r)) |v'(r)|^2+\frac{n-1}{r^2}|v(r)|^2 \right)r^{n-1}\,dr
\end{equation}
with the boundary condition
\begin{equation}\label{variationalBC}
v(\rho)=\left( \meanint_{S^{n-1}} |u(\rho\,\theta)|^2\,d\theta\right)^{1/2}.
\end{equation}
\end{subequations}
If we let $v_*$ denote the minimizer for \eqref{variational}, then we can express it as a scalar multiple of $Z(r)$ by
\begin{equation}\label{v*=}
v_*(r)=\frac{Z(r)}{Z(\rho)}\left( \meanint_{S^{n-1}} |u(\rho\,\theta)|^2\,d\theta\right)^{1/2},
\end{equation}
and the minimum value of $J$ is
\begin{equation}\label{J(v*)=}
J(v_*)=(1+g(\rho))\rho^{n-1}\frac{Z'(\rho)}{Z(\rho)}\meanint_{S^{n-1}} |u(\rho\,\theta)|^2\,d\theta.
\end{equation}

Now if we fix $\theta\in S^{n-1}$ and evaluate $J$ on $v(r)=u(r\theta)$, we obtain
\[
\int_0^\rho \left( (1+g(r))\left|\frac{\partial v}{\partial r}\right|^2+\frac{n-1}{r^2}\left| v\right|^2\right)\,r^{n-1}dr \geq J(v_*).
\]
If we integrate both sides over $\theta\in S^{n-1}$ and use \eqref{estimate} and \eqref{J(v*)=}, we obtain
\begin{equation}\label{geq_Z'/Z}
\frac{1}{2} \frac{d}{d\rho}\|u(\rho)\|^2\geq \frac{Z'(\rho)}{Z(\rho)} \|u(\rho)\|^2 \quad\hbox{where}\ \|u(\rho)\|^2:=\int_{S^{n-1}} |u(\rho\,\theta)|^2\,d\theta.
\end{equation}
But \eqref{geq_Z'/Z} implies
\[
\frac{d}{d\rho}\log\frac{\|u(\rho)\|}{Z(\rho)}\geq 0,
\]
which implies that $\|u(\rho)\|/Z(\rho)$ is increasing, so \eqref{eq:comparison} follows. $\Box$

\subsection{Conditions for Loss of Lipschitz Regularity}
Propositions \ref{pr:1} and \ref{pr:2} provide conditions under which Lipschitz regularity holds. Now let us consider conditions under which a weak solution exists which is {\it not} Lipschitz continuous.
\begin{proposition}\label{pr:NotLip} 
Let the coefficients $a_{ij}$ be given by  \eqref{GS}, where $g(r)$ is square-Dini continuous at $0$ and 
\begin{equation}\label{int-g(r)/r=->infty}
\limsup_{r\to 0}\int_{r}^{1/2} \frac{g(\rho)}{\rho}\,d\rho =+\infty.
\end{equation}
Also assume $g(r)$ has finite total variation, i.e.\ $\int_0^{1/2}|g'(r)|\,dr<\infty$.
Then there exists a weak solution $u$ of \eqref{Lu=0} in $B_\rho$ that is not Lipschitz continuous at $x=0$.
If we also have $g(r)>0$ for $r<\e$, then we do not require $g(r)$ to be square-Dini continuous at $r=0$ for the conclusion to hold.
\end{proposition}

 \noindent  
{\bf Proof.} We want to find a weak solution $u\in H^{1,2}(B_\rho)$ of \eqref{Lu=0} that is bounded but not Lipschitz continuous. Without loss of generality, we may assume $\rho=1/2$ and $u(0)=0$. We shall look for a solution 
in the form $u(r,\theta)=v(r)\,\Theta_1(\theta)$, where $v(0)=0$ and $\Theta_1$ is an eigenfunction for the first eigenvalue $\lambda=n-1$ for the Laplace-Beltrami operator on the unit sphere. We find that $v(r)$ should satisfy the ordinary differential equation \eqref{ode:v}.
With a change of independent variable to $t=-\log r$, the equation for $v$ can be written as
\begin{equation}\label{ode:v(t)}
\frac{d}{dt}\left( (1+{\rm g})\frac{dv}{dt}\right)-(n-2)(1+{\rm g})\frac{dv}{dt}-(n-1)v=0, \quad\hbox{for $0<t<\infty$,}
\end{equation}
where ${\rm g}(t)=g(e^{-t})$, and this can be written as a 1st-order system
\begin{equation}\label{dynsys-v,w}
\frac{d}{dt}\begin{bmatrix} v \\ w \end{bmatrix} = \begin{bmatrix} 0 & (1+{\rm g})^{-1} \\ n-1 & n-2 \end{bmatrix}\begin{bmatrix} v \\ w \end{bmatrix} \quad\hbox{where} \ w=(1+{\rm g})\frac{dv}{dt}.
\end{equation}
Let $M_{\rm g}$ denote the $2\times 2$ matrix on the right hand side of \eqref{dynsys-v,w}.
For ${\rm g}=0$,  $M_0$  has eigenvalues $\lambda_1=-1$ and $\lambda_2=n-1$ (corresponding to the solutions $r$ and $r^{1-n}$ of \eqref{ode:v} when $ g=0$). If we let $E$ denote the matrix of eigenvectors for  $M_0$,
then we can rewrite \eqref{dynsys-v,w} as
\begin{subequations}\label{eq:dynsys-V}
\begin{equation}
\frac{d}{dt}V =
(A_0 + B({\rm g}) )\,V, \quad\hbox{where} \ \begin{bmatrix} v \\ w \end{bmatrix} = E\,V
\end{equation}
and
\begin{equation}
A_0= \begin{bmatrix} -1 & 0 \\ 0 & n-1 \end{bmatrix}, \quad
B({\rm g})= \frac{\rm g}{n(1+{\rm g})}\begin{bmatrix} n-1 &  -(n-1)^2\\  1 & -(n-1)
\end{bmatrix}. \end{equation}
\end{subequations}
The eigenvalues of $A_0+B({\rm g})$ may be denoted  $\Lambda_1({\rm g})$ and $\Lambda_2({\rm g})$, where
$\Lambda_1({\rm g})\to -1$ and $\Lambda_2({\rm g})\to n-1$ as ${\rm g}\to 0$. We can calculate $\Lambda':=d\Lambda/d{\rm g}$ and evaluate $\Lambda'(0)$ to find
\begin{equation}\label{lambda'(0)=}
\Lambda_1'(0)=\frac{n-1}{n}.
\end{equation}
Consequently
\begin{equation}
\Lambda_1({\rm g})=-1+\frac{n-1}{n}{\rm g}+O(|{\rm g}|^2) \quad\hbox{as}\ {\rm g}\to 0.
\end{equation}
Similarly, we can use $\Lambda(0)=n-1$ and obtain
\begin{equation}
\Lambda_2({\rm g})=n-1+\frac{n-1}{n}{\rm g}+O(|{\rm g}|^2) \quad\hbox{as}\ {\rm g}\to 0.
\end{equation}

At this point we have not used the fact that $g(r)$ is square-Dini continuous at $r=0$, so we do not assert that
$\int_0^\infty |{\rm g}(t)|^2\,dt<\infty$. Nevertheless, because we have assumed that $g(r)$ has finite total variation for
$0<r<1/2$, we know that the matrix function $B(t):=B({\rm g}(t))$ has finite total variation for $t>\log 2$, so we can
use the asymptotic theory of linear systems of ordinary differential equations.
Denoting the eigenvalues as $\lambda_1(t):=\Lambda_1({\rm g}(t))$ and $\lambda_2(t):=\Lambda_2({\rm g}(t))$, we have
\begin{equation}\label{lambda1-asym}
\lambda_1(t)=-1+\frac{n-1}{n}\,{\rm g}(t) + O(|{\rm g}(t)|^2) \quad\hbox{as} \ t\to\infty.
\end{equation}
We can apply Theorem 8.1 in \cite{CL} (see also Theorem 16.3.1 in \cite{KM}) to conclude that there exists a solution $V_1(t)$ with the asymptotic behavior
\begin{equation}\label{V1-asym1}
\lim_{t\to\infty} V_1(t)\,\exp\left[-\int_0^t \lambda_1(\tau)\,d\tau\right] =V_0\, \quad
\hbox{as}\ t\to\infty,
\end{equation}
where $V_0=(1,0)^t$ is the eigenvector for the eigenvalue $\Lambda_1(0)=-1$.

Now we shall invoke the sqaure-Dini assumption $\int_0^\infty |{\rm g}(t)|^2\,dt<\infty$. This mean that \eqref{eq:dynsys-V} can be written
\begin{subequations}\label{eq:dynsys-V2}
\begin{equation}\label{eq:dynsys-V2a}
\frac{d}{dt}V =
(A_0 + {\rm g}(t)\, B_0) \,V + C( {\rm g}(t))\,V,
\end{equation}
where 
\begin{equation}
A_0= \begin{bmatrix} -1 & 0 \\ 0 & n-1 \end{bmatrix}, \quad
B_0:= \frac{1}{n}\begin{bmatrix} n-1 &  -(n-1)^2\\  1 & -(n-1)
\end{bmatrix}, \quad\hbox{and}\quad
 \int_0^\infty |C( {\rm g}(t))|\,dt<\infty.
\end{equation}
\end{subequations}
Now  Theorem 8.1 in \cite{CL} shows that
\eqref{V1-asym1} can be sharpened somewhat to
\begin{equation}\label{V1-asym2}
V_1(t)=V_0\,\exp\left[-t+\frac{n-1}{n}\int_0^t{\rm g}(\tau)\,d\tau\right](1+o(1)) \quad
\hbox{as}\ t\to\infty.
\end{equation}
If we transform back to our original dependent variables $v,w$ using $(v,w)^t=E^{-1}V$ and return to the original independent variable $r=e^{-t}$, then we find there is a solution $v(r)$ of \eqref{ode:v} with the asymptotic behavior
\begin{equation}\label{v,w-asym}
v(r)= c\,r\exp\left[ \frac{n-1}{n}\int_r^{1/2} \frac{g(\rho)}{\rho}\,d\rho\right](1+o(1)) \quad\hbox{as}\ r\to 0.
\end{equation}
But the condition \eqref{int-g(r)/r=->infty} shows that this solution of \eqref{ode:v} is not Lipschitz at $r=0$.
In other words, we have found a bounded solution $u(r,\theta)=v(r)\Theta_1(\theta)$ of \eqref{Lu=0} that is not Lipschitz at $r=0$.

If we do not assume that $g(r)$ is square-Dini at $r=0$, but we know $g(r)>0$ for $r<\e$, then we still can conclude
\[
\frac{n-1}{n}\int_0^t \left({\rm g}(\tau)+O(|{\rm g}(\tau)|^2)\right)\,d\tau\to +\infty \quad\hbox{as}\ t\to\infty,
\]
which shows that our solution $v(r)$ is not Lipschitz continuous at $r=0$.
$\Box$

\begin{remark}\label{re:NASC} 
 Propositions \ref{pr:2} and \ref{pr:NotLip} both depend on the asymptotic behavior of the positive finite energy solution $Z(r)$ of \eqref{ode:v} in such a way that $Z(r)\leq c\,r$ becomes a necessary and sufficient condition for every weak solution of \eqref{Lu=0} to be Lipschitz continuous at $x=0$.
\end{remark}

\begin{remark}\label{re:2} 
As we shall see in the following Examples, it is possible to find functions $g(r)$ that do not have Dini mean oscillation, yet satisfy the conditions of Proposition \ref{pr:1} or \ref{pr:2} that show all weak solutions are $C^1$ in a neighborhood of $x=0$.
Proposition \ref{pr:NotLip}, on the other hand, can be used to find examples of coefficients $a_{ij}$ that do not satisfy the Dini mean oscillation condition: if $g(r)$  satisfies the square-Dini condition and \eqref{int-g(r)/r=->infty}, then \eqref{Lu=0} has weak solutions $u\in H^{1,2}(\Omega)$ that are not Lipschitz continuous, so by the results of \cite{DK} the coefficients \eqref{GS} cannot have Dini mean oscillation.
\end{remark}

\subsection{Examples}
In this section we consider several specific functions $g(r)$ and what our preceding results tell us about Lipschitz regularity at $x=0$ or in a neighborhood of $x=0$. We frequently find that Proposition \ref{pr:2} enables us to extend results beyond what Proposition \ref{pr:1} is able to provide due to the condition that $g(r)$ be square-Dini continuous at $r=0$. But both Propositions  \ref{pr:1} and  \ref{pr:2} apply in some cases that do not have Dini mean oscillation.

\medskip\noindent
{\bf Example 1.} 
Let us consider \eqref{GS} with
\begin{equation}\label{g(r)=}
g(r)=|\log r|^{-\gamma} \quad \hbox{for $\gamma>0$.}
\end{equation}  
For $\gamma>1$, $g(r)$ is Dini continuous at $r=0$, so we know weak solutions of \eqref{Lu=0} are Lipschitz continuous (in fact, $C^1$).
For $0<\gamma\leq 1$, $g(r)$ is not Dini continuous: in fact
\[
\int_0^{1/2}\frac{g(r)}{r}\,dr=\int_0^{1/2} \frac{|\log r|^{-\gamma}}{r}\,dr=+\infty,
\]
so condition \eqref{int-g(r)/r=->infty}  holds. Moreover, $g(r)$ has finite total variation:
\[
\int_0^{1/2} |g'(r)|\,dr=\gamma\int_0^{1/2} \frac{|\log r|^{-\gamma-1}}{r}\,dr <\infty.
\]
Now for $1/2<\gamma\leq 1$, $g(r)$ is square-Dini continuous, so we can use Proposition \ref{pr:NotLip}
to conclude there exists a weak solution of  \eqref{Lu=0} that is not Lipschitz continuous at $x=0$. 
In fact, since $g(r)>0$ for $0<r<1/2$ we can use Proposition \ref{pr:NotLip} to extend this conclusion to $0<\gamma\leq 1$ .
 {\bf Summary for \eqref{g(r)=}}: {\it For $\gamma>1$, all weak solutions of \eqref{Lu=0} are Lipschitz continuous; but for $0<\gamma\leq 1$, there exists a weak solution of \eqref{Lu=0} that is not Lipschitz continuous at $x=0$.}

Now let us  consider  \eqref{GS} with
\begin{equation}\label{g(r)=-}
g(r)=-|\log r|^{-\gamma} \quad \hbox{for $\gamma>0$.}
\end{equation}  
For $\gamma>1$, as with \eqref{g(r)=}, $g(r)$ is  Dini continuous at $r=0$, so we know weak solutions are Lipschitz continuous; but for all $\gamma>0$ we have $|rg'(r)|\leq c$ and
\[
-\infty\leq \int_0^{1/2}\frac{g(r)}{r}\,dr=-\int_0^{1/2}\frac{|\log r|^{-\gamma}}{r} <0.
\]
For $1/2<\gamma\leq 1$, $g(r)$ is square-Dini continuous and  the condition \eqref{cond:g(r)/r} is trivially satisfied with $K=0$, so we could use Proposition \ref{pr:1} to conclude that all weak solutions of \eqref{Lu=0} are Lipschitz continuous. However, we can instead use Proposition \ref{pr:2} to handle all values $0<\gamma\leq 1$. We need to find the solution $Z(r)$ of \eqref{ode:v} with the desired properties. 
From \eqref{lambda1-asym} and \eqref{V1-asym1} we see that there is a solution $V(t)$ of \eqref{eq:dynsys-V} such that $|V(t)|= o(e^{-t})$, which means we have our solution $v(t)=Z(r)$ satisfying $|Z(r)|\leq c\,r$. If we rescale so that $Z(1)=1$, then the non-oscillatory character of \eqref{ode:v} shows $Z(r)$ remains positive. We can also confirm that $Z$ satisfies the finite energy condition \eqref{finite_energy}, so we can apply Proposition \ref{pr:2} to reach our conclusion. {\bf Summary for \eqref{g(r)=-}}: {\it For $\gamma>0$, all weak solutions of \eqref{Lu=0} are Lipschitz continuous.}

\smallskip\noindent{\bf Note:}
For $0<\gamma\leq 1$, these results for \eqref{g(r)=} and \eqref{g(r)=-}  {cannot} be obtained using Dini mean oscillation for $a_{ij}$ because that condition is satisfied  only for $\gamma>1$: this can be confirmed by direct calculation at $x=0$ (see the Appendix), but for \eqref{g(r)=} it also follows from Remark \ref{re:2}.

 \bigskip\noindent
{\bf Example 2.} 
Let us consider an example where $g(r)$ oscillates as $r\to 0$:
\begin{equation}\label{g(r)=sin}
g(r)=\frac{\sin(|\log r|)}{|\log r|^\beta}, \quad\hbox{where $\beta>0$.}
\end{equation}  
The condition $|rg'(r)|\leq c$ in Propositions \ref{pr:1} and \ref{pr:2} is satisfied. 
Also note that
\begin{equation}\label{Ex2:int(g/r)}
\int_0^{e^{-1}} \frac{g(r)}{r}\,dr
=\int_1^\infty \frac{\sin\tau}{\tau^{\beta}}\,d\tau \quad\hbox{converges}
\end{equation} 
since $\tau^{-\beta}$
is strictly decreasing and $\sin\tau$ alternates sign on intervals $[k\pi,(k+1)\pi]$ for $k=0,1,\dots$. In fact, for any interval $[t_1,t_2]\subset (2\pi,\infty)$ we have 
\begin{equation}\label{Ex2:stability}
\int_{t_1}^{t_2}\frac{\sin \tau}{\tau^{\beta}}\,d\tau \leq \int_{2\pi}^{3\pi} \frac{\sin \tau}{\tau^{\beta}}\,d\tau <\infty.
\end{equation}  
For $\beta>1/2$, this $g(r)$ is square-Dini continuous and \eqref{Ex2:stability}
 confirms condition \eqref{cond:g(r)/r}, so by Proposition \ref{pr:1}  we conclude that all weak solutions of \eqref{Lu=0} are Lipschitz continuous. To obtain this conclusion for all $\beta>0$, let us appeal to Proposition \ref{pr:2}. Using asymptotic analysis as above with \eqref{Ex2:int(g/r)}, we can find a finite energy solution $Z(t)$ of \eqref{ode:v}
satisfying $0<Z(r)\leq c\,r$ for $0<r<1$, so Proposition \ref{pr:2} shows that all weak solutions are Lipschitz continuous.

\smallskip\noindent{\bf Note:} These results for  $0<\beta\leq 1$ in \eqref{g(r)=sin}
 { cannot} be obtained using Dini mean oscillation since one can show by direct calculation (see the Appendix) that the coefficients $a_{ij}$ have Dini mean oscillation only for $\beta>1$.

\bigskip\noindent
{\bf Example 3.} 
Consider the following function:
\begin{equation}\label{Ex3:g(r)=}
g(r)=\frac{-C_1\,\sin(|\log r|)-C_2\cos(|\log r|)}{A+\sin(|\log r|) -\cos(|\log r|)}
\end{equation}
where
\[
C_1=\frac{(n-1)^2}{(n-1)^2+1}+1, \quad C_2=\frac{n-1}{(n-1)^2+1}-1,
\]
and $A>1$  is chosen sufficiently large. This complicated expression for $g(r)$ arises because we actually want the following function $Z(r)$ to be a solution of \eqref{ode:v}:
\begin{equation}
Z(r)=r\,\left(A+\sin(|\log r|)\right).
\end{equation}
Using $Z(r)$ in \eqref{ode:v}, we can work backwards to find $g(r)$ and confirm \eqref{Ex3:g(r)=}.
This means that, if we use this $g(r)$ in our coefficients $a_{ij}$ as in \eqref{GS}, then by Proposiiton \ref{pr:2} any weak solution $u\in H^{1,2}(B_{1/2}(0))$ of \eqref{Lu=0} must be Lipschitz continuous in $B_{1/2}(0)$.
 This is in spite of the fact that the coefficients
 are not square-Dini continuous  or Dini mean oscillation (see Appendix);  in fact, $g$ does not even vanish as $r\to 0$, and we claim that
 \begin{equation}\label{Ex3:int=+infty}
 \int_0^{1/2} \frac{g(r)}{r}\,dr=+\infty.
 \end{equation}
This is interesting since the condition \eqref{Ex3:int=+infty} is associated with the existence of a weak solution that is not Lipschitz continuous (cf.\ Proposition \ref{pr:NotLip}). However, since $g(r)$ does not even vanish as $r\to 0$, any intuition from the  asymptotic analysis used in the proof of Proposition \ref{pr:NotLip}
 is not relevant for this example.

To confirm \eqref{Ex3:int=+infty}, let us compute
\[
\int_0^{1} \frac{g(r)}{r}\,dr=\int_{0}^\infty {\rm g}(t)\,dt
=\int_{0}^\infty \frac{-C_1\,\sin t-C_2\cos t}{A+\sin t -\cos t}\,dt
\]
By periodicity, it suffices to compute
\[
\begin{aligned}
\int_0^{2\pi} \frac{-C_1\,\sin t-C_2\cos t}{A+\sin t -\cos t}\,dt & =
\frac{1}{A}\int_0^{2\pi} \frac{-C_1\,\sin t-C_2\cos t}{1+\frac{1}{A}(\sin t -\cos t)}\,dt \\
& = A^{-2}\int_0^{2\pi} (\sin t-\cos t) \left[C_1\sin t+C_2\cos t\right]\,dt + O(A^{-3})
\end{aligned},
\]
provided $A$ is sufficiently large. By orthogonality, the cross terms vanish and we are left with computing
\[
\int_0^{2\pi}  \left[C_1\sin^2 t-C_2\cos^2 t\right]\,dt
=\int_0^{2\pi} \frac{C_1-C_2}{2}\,dt,
\]
where we have used  half-angle trigonometric identities for $\sin^2$ and $\cos^2$ as well as $\int_0^{2\pi}\cos 2t\, dt=0
=\int_0^{2\pi}\sin 2t\,dt$. But it is clear that $C_1-C_2>0$, so we obtain
\[
\int_{0}^\infty \frac{-C_1\,\sin t-C_2\cos t}{A+\sin t -\cos t}\,dt = +\infty.
\]

 \appendix
\section{Computations of  mean oscillation at $x=0$}\label{sec:A}
All of the functions that we consider are smooth except at $x=0$, so we will restrict our attention
to the mean oscillation at $x=0$, i.e.\ we take the sup over $B_r=\{y:|y|<r\}$ in \eqref{Dini2}. However, if we want to show  for certain examples that Dini mean oscillation over a domain fails, it suffices to show that it fails at $x=0$.

First let us confirm that $a(x)=|\log r|^{-\gamma}$ has Dini mean oscillation at $x=0$ for all $\gamma>0$.
Since $a(x)=a(r)$ is radial, we easily compute
\begin{equation}\label{MV(a)}
\begin{aligned}
\widetilde a(r):=\widetilde a_{0,r}&=\meanint_{B_r}a(x)\,dx=\frac{n}{r^n}\int_0^r |\log \rho|^{-\gamma}\,\rho^{n-1}\,d\rho \\
&=|\log r|^{-\gamma} -\gamma\,|\log r|^{-\gamma-1}+\frac{\gamma(\gamma+1) n}{r^n}\int_0^r |\log \rho|^{-\gamma-2}\,\rho^{n-1}\,d\rho.
\end{aligned}
\end{equation}
So
$
a(x)-\widetilde a(r)=\gamma\,|\log r|^{-\gamma-1}+O(|\log r|^{-\gamma-2})
$
and hence the mean oscillation is
\begin{equation}
\omega_a(r)=\meanint_{B_r}|a(x)-\widetilde a(r)|\,dx= \gamma\,|\log r|^{-\gamma-1}+O(|\log r|^{-\gamma-2}).
\end{equation}
Since $\om_a(r)$ satisfies \eqref{Dini1}, this confirms that $a(x)$ has Dini mean oscillation at $x=0$ for all $\gamma>0$. 

In the above calculation, there was no angular oscillation, and this contributed to $a$ satisfying the Dini mean oscillation condition for $\gamma>0$. The situation is quite different for functions with angular dependence. 
In order to estimate the mean oscillation for the coefficients $a_{ij}$ in \eqref{GS}, let us introduce the matrix
\begin{equation}\label{def:Theta}
\Theta = \left( \theta_i\theta_j \right)_{i,j=1,\dots,n},
\end{equation}
and denote its mean value over the unit sphere $S^{n-1}$ by $\widetilde \Theta$. Then a calculation shows
\begin{equation}\label{TildeTheta}
\widetilde \Theta = \frac{1}{n} \left( \delta_{ij} \right)_{i,j=1,\dots,n}.
\end{equation}
Now let us multiply the matrix $\Theta$ by a scalar function $g(r)$ that is smooth for $0<r<1/2$. If we denote by $\widetilde{g\Theta}(r)$ the mean value of $g(r)\Theta$ over the ball $B_r(0)$, then we have 
\[
\widetilde{g\Theta}(r)=\widetilde g(r) \,\widetilde\Theta, \quad\hbox{where}\ \widetilde g(r)=\frac{n}{r^n}\int_0^r g(\rho)\rho^{n-1}\,d\rho.
\]
Now we may estimate the matrix norm of $g\,\Theta-\widetilde{g\,\Theta}$ from below by 
\begin{equation}\label{est:gTheta}
\pmb{|} g(r)\Theta-\widetilde{g\,\Theta}(r) \pmb{|} \geq 
  |g(r)|\, \pmb{|} \Theta-\widetilde{\Theta} \pmb{|} - |g(r)-\widetilde g(r)|\, \pmb{|} \widetilde{\Theta} \pmb{|}.
\end{equation}
If we can show $|g(r)-\widetilde g(r)|=o(|g(r)|)$ as $r\to 0$, then we know that the matrix $g(r)\Theta$ has Dini mean oscillation at $x=0$ if and only if the function $|g(r)|$ itself satisfies the Dini condition at $r=0$.

For example, consider the coefficients in Example 1, which we can write in matrix form as
\begin{equation}
A(x)=I+|\log r|^{-\gamma}\Theta, \quad\hbox{where $I$ denotes the identity matrix}.
\end{equation}
The mean value of $I$ is itself, so we have
\[
A(x)-\widetilde A(r)=g(r)\,\Theta -\widetilde{g}(r)\,\widetilde\Theta \quad\hbox{where}\ g(r)=|\log r|^{-\gamma}.
\]
As we found in \eqref{MV(a)}, $|g(r)-\widetilde g(r)|=o(|g(r)|)$ as $r\to 0$, so by \eqref{est:gTheta}, $A(x)$ has Dini mean oscillation at $x=0$ if and only if $g(r)=|\log r|^{-\gamma}$ satisfies the Dini condition. So Example 1 has Dini mean oscillation at $x=0$ for
$\gamma>1$, but not for $0<\gamma\leq 1$.

Next let us consider the coefficients in Example 2. Again we need only concern ourselves with the matrix
\[
g(r)\Theta \quad\hbox{where}\ g(r)=\frac{\sin(|\log r|)}{|\log r|^\beta}.
\]
For $0<r<1$ we use integration by parts and a geometric series to evaluate
\[
\begin{aligned}
\widetilde g(r)&=\frac{n}{r^{n-1}}\int_0^r \frac{\sin(-\log\rho)}{(-\log \rho)^\beta}\rho^{n-1}\,d\rho
=n e^{nt}\int_t^\infty \frac{\sin\tau}{\tau^\beta}e^{-n\tau}d\tau \\
&=\frac{\sin t}{t^\beta}+e^{nt}\int_t^\infty \frac{\cos\tau}{\tau^\beta} e^{-n\tau}d\tau 
=\frac{\sin t}{t^\beta}+\frac{\cos t}{nt^\beta}-\frac{1}{n}e^{nt}\int_t^\infty \frac{\sin\tau}{\tau^\beta}e^{-n\tau}d\tau\\
&=\frac{n}{n^2+1}\left(\frac{n\sin t+\cos t}{t^\beta}\right).
\end{aligned}
\]
Consequently,
\begin{equation}\label{Ex2:g-tilde(g)}
g(r)-\widetilde g(r) = \frac{1}{n^2+1}\,\frac{\sin(|\log r|)-n\cos(|\log r|)}{|\log r|^\beta}.
\end{equation}  
Now we can estimate
\[
\meanint_{S^{n-1}} \pmb{|} g(r)\Theta-\widetilde g(r) \widetilde \Theta \pmb{|}\,ds
\geq \left|\, \meanint_{S^{n-1}} ( g(r)\Theta-\widetilde g(r) \widetilde \Theta )\,ds\, \right|
=|g(r)-\widetilde g(r)|\, \pmb{|}\widetilde\Theta \pmb{|}.
\]
But if $0<\beta\leq 1$, then
\[
\int_0^{1/2} \frac{|g(r)-\widetilde g(r)|}{r}\,dr=\frac{1}{n^2+1}\int_{\log 2}^\infty \frac{|\sin t-n\cos t|}{t^\beta}\,dt=\infty
\]
 since the numerator is periodic. This shows that  the coefficients in Example 2 have Dini mean oscillation at $x=0$
only for $\beta>1$.

Finally, we consider the coefficients in Example 3. Recall 
\[
g(r)=\frac{-C_1\sin (|\log r|)-C_2\cos(|\log r|)}{A+\sin(|\log r|) -\cos(|\log r|)}
= A^{-1} (-C_1\sin (|\log r|)-C_2\cos(|\log r|))+O(A^{-2})
\]
for large $A$.
We can use this to approximate $\widetilde g(r)$, the mean value over the ball $B_r$, and then express the answer in terms of $t=-\log r$:
\[
\widetilde g(e^{-t})= -\frac{n}{A(n^2+1)} (C_1(n\sin t-\cos t)+C_2(n\cos t-\sin t))+O(A^{-2}).
\]
Hence
\[
{ g}(e^{-t})-\widetilde{g}(e^{-t})=
-\frac{1}{A(n^2+1)}\left[ (C_1+nC_2)\sin t+(C_2+nC_1)\cos t\right]+O(A^{-2}).
\]
But
\[
\int_{\log r}^\infty |(C_1+nC_2)\sin t+(C_2+nC_1)\cos t|\,dt
=+\infty,
\]
since the integrand is periodic, so $|g(r)-\widetilde g(r)|$ does not satisfy the Dini condition. Hence Example 3 does not have Dini mean oscillation.

\bigskip\noindent
{\bf Acknowledgement}: This paper has been supported by the RUDN University Strategic
Academic Leadership Program.



\begin{thebibliography}{99}   

\bibitem{ADN} S. Agmon, A. Douglis and L. Nirenberg, \textit{Estimates near the boundary for solutions of elliptic partial dif- 
ferential equations satisfying general boundary conditions I}, Comm. Pure Appl. Math., \textbf{12} (1959), 623-727.

\bibitem{A} Yu. A. Alkhutov,  \textit{Smoothness and limiting properties of solutions of a second-order parabolic equation}, 
Matematicheskie Zametki, \textbf{50}, 150-152.

\bibitem{B} C. Burch, \textit{The Dini condition and regularity of weak solutions of elliptic equations}, J. Diff. Eq. \textbf{30} (1978), 308-323.

\bibitem{CFL} F. Chiarenza, M. Frasca and P. Longo,  \textit{Solvability of the Dirichlet problem for nondivergence elliptic equa- tions with VMO coefficients}, Trans. Amer. Math. Soc.,  \textbf{336} (1993), 841–853.

\bibitem{CM}  A. Cianchi, V. Maz’ya, \textit{Global Lipschitz regularity for a class of quasilinear elliptic }\, Commun. PDE \textbf{36},  (2011), 100–133.

\bibitem{CM1}  A. Cianchi, V. Maz’ya, \textit{Gradient regularity via rearrangements for p-Laplacian type elliptic boundary value problems}, Eur. Math. Soc. (JEMS), (2014), \textbf{16}, p.571-595

\bibitem{CM2}  A. Cianchi, V. Maz’ya, \textit{Global boundedness of the gradient for a class of nonlinear elliptic systems},
Arch. Ration. Mech. Anal., \textbf{212} (2014), p.129-177.

\bibitem{CM3}  A. Cianchi, V. Maz’ya, \textit{Global gradient estimates in elliptic problems under minimal data and domain regularity}, Commun. Pure Appl. Anal. \textbf{14} (2015), no. 1, 285–311.

\bibitem{CL} E. Coddington, N. Levinson, \textit{Theory of Ordinary Differential Equations}, McGraw-Hill, New York, 1955.

\bibitem{DM} C. De Filippis, G. Mingione, \textit{Lipschitz bounds and nonautonomous integrals}, Arch. Rat. Mech. Anal, doi: 10.1007/s00205-021-01698-5

\bibitem{DG} E. De Giorgi, \textit{Sulla differenziabilit\`a e l'analiticit\`a delle estremali degil integrali multipli regolari}, Mem. Accad. Sci. Torino Cl. Sci. Fis. Mat. Natur. (3) \textbf{3} (1957), 25-43.

\bibitem{DF} G. Di Fazio, \textit{Estimates for divergence form elliptic equations with discontinuous coefficients}, Boll. Un. 
Mat. Ital A(7), \textbf{10} (1996), 409-420. 

\bibitem{DEK1} H. Dong, L. Escauriaza, S. Kim, \textit{On $C^1$, $C^2$, and weak type $(1,1)$ estimates for linear elliptic operators, Part II}, Math. Ann. \textbf{370} (2018), no. 1-2, 447-489.

\bibitem{DEK2} H. Dong, L. Escauriaza, S. Kim, \textit{On $C^{1/2,1}$, $C^{1,2}$, and $C^{0,0}$ estimates for linear parabolic  operators}, J. Evol. Equ. (2021) \ https://doi.org/10.1007/s00028-021-00729-8

\bibitem{DK} H. Dong, S. Kim, \textit{On $C^1$, $C^2$, and weak type $(1,1)$ estimates for linear elliptic operators}, Comm. Partial Differential Equations, \textbf{42} (2017), 417-435.

\bibitem{DLK} H. Dong, J. Lee, S. Kim, \textit{On conormal and oblique derivative problem for elliptic equations with Dini mean oscillation coefficients}, Indiana Univ. Math. J. \textbf{69} (2020), no. 6, 1815–1853. 


\bibitem{GS} D. Gilbarg, J. Serrin,  \textit{On isolated singularities of solutions of second-order elliptic equations,} J.\ Analyse Math.\ \textbf{ 4} (1955/56), 309-340.


\bibitem{HW} P. Hartman, A. Wintner, \textit{On uniform Dini conditions in the theory of linear partial differential equations of elliptic type,} American J. Math. \textbf{77} (1955), 329-354.

\bibitem{I} A.M. ll'in, \textit{On parabolic equations whose coefficients satisfy Dini's condition}, Mat. Zametki, \textbf{I}  (1967), 71-80 = Math. Notes \textbf{1} (1967), 46- 51.

\bibitem{JMV} T. Jin, V.G. Maz'ya, J. Van Schaftingen, \textit{Patholgical solutions to elliptic problems in divergence form with continuous coefficients}, Comptes rendus. Math\'ematique, \textbf{47} (2009), p.773-778.

\bibitem{KM} V.A. Kozlov, V.G. Maz'ya, \textit{Differential Equations with Operator Coefficients}, Springer-Verlag, New York, 1999.

\bibitem{KM2}  V.A. Kozlov, V.G. Maz’ya, \textit{Asymptotic formula for solutions to the Dirichlet problem of elliptic equations with discontinuous coefficients near the boundary}, Ann. Scuola Norm. Sup. Pisa Cl. Sci. (5), \textbf{2} (2003), 551-600.

\bibitem{KuM} T. Kuusi, G. Mingione, \textit{Universal potential estimates}, J. Funct. Anal. \textbf{262} (2012), 4205-4269.

\bibitem{KuM2}  T. Kuusi, G. Mingione, \textit{A nonlinear Stein theorem}, Calc. Var. \textbf{51} (2014), 45–86.

\bibitem{L} E. M. Landis, \textit{Second order equations of elliptic and parabolic type}, Translations of Mathematical Monographs, 171. American Mathematical Society, Providence, RI, 1998.


\bibitem{ME} M. I. Matıchuk,  S. D. Eıdel'man,  \textit{On parabolic systems with coefficients satisfying Dini’s condition}. Dokl. Akad. Nauk SSSR  \textbf{165} (1965), 482-485.

\bibitem{M} V.G. Maz'ya, \textit{On weak solutions of the Dirichlet and Neumann problems}, Trans.\ Moscow Math. Soc. \textbf{20} (1969), 135-172.

\bibitem{MM} V.G. Maz'ya, R.McOwen,  \textit{Differentiablilty of solutions to second-order elliptic equations via dynamical systems}, J. Differential Equations, \textbf{250} (2010), 1137-1168.

\bibitem{MM1} V.G. Maz'ya, R.McOwen, \textit{Second-order differentiability for solutions of elliptic equations in the plane},  J. Mathematical Sciences, 191 (2013), 243-253.

\bibitem{MM2} V.G. Maz'ya, R.McOwen, \textit{Differentiability of solutions to the Neumann problem with low-regularity data via dynamical systems}, {\it Operator Theory: Advances and Applications}, 261 (2017), 343-385.

 \bibitem{MMS} V.G. Maz'ya, M. Mitrea, T. Shaposhnikova,  \textit{The Dirichlet problem in Lipschitz domains for higher order elliptic systems with rough coefficients}, J. Anal. Math. \textbf{110} (2010), 167–239. 
 
 \bibitem{MS} V.G. Maz'ya,  T. Shaposhnikova,   \textit{Theory of Sobolev Multipliers}, Springer, 2009.
 
 \bibitem{Mc} R. McOwen, \textit{On elliptic operators in nondivergence and double divergence form}, Operator Theory: Advances and Applications, 193 (2009), 159-169.

\bibitem{Mo} J. Moser,  \textit{A new proof of De Giorgi’s theorem concerning the regularity problem for elliptic differential equations}, Comm. Pure Appl. Math. \textbf{13} (1960), 457–468.

\bibitem{N} J. Nash, \textit{Continuity of solutions of parabolic and elliptic equations}, Amer. J. Math. \textbf{80} (1958), 931-954.

\bibitem{Se} J. Serrin, \textit{ Pathological solutions of elliptic differential equations}, Ann. Scuola Norm. Sup. Pisa (3) {\bf 18} (1964), 385–387.


 \bibitem{T} M. Taylor,  \textit{Tools for PDE}, Mathematical Surveys and Monographs {\bf 81}, AMS, 2000.
 
 
\end{thebibliography}
\end{document}